\documentclass[10pt,twoside]{amsart}

\usepackage{amssymb}
\usepackage{amsxtra}
\usepackage{graphicx}
\usepackage{dcpic,pictexwd}

\theoremstyle{plain}

\newtheorem{thm}{Theorem}
\newtheorem{theorem}{Theorem}[section]
\newtheorem{prop}[theorem]{Proposition}

\newtheorem{defn}[theorem]{Definition}

\newtheorem{lemma}{Lemma}

\newcommand{\al}{\alpha}
\newcommand{\be}{\beta}
\newcommand{\ga}{\gamma}

\newcommand{\Q}{\mathbb{Q}}

\newcommand{\Z}{\mathbb{Z}}

\newcommand\ra{\rightarrow}
\newcommand\lra{\longrightarrow}

\newcommand{\Si}{\Sigma}

\newcommand{\ben}{\begin{enumerate}}
\newcommand{\een}{\end{enumerate}}
\newcommand{\disp}{\displaystyle}

\newcommand{\Oz}{P.\ Ozsv{\'a}th\,}
\newcommand{\Sz}{Z.\ Szab{\'o}\,}

\begin{document}
\title{Some bounds for the knot Floer $\tau$-invariant of satellite knots}
\author{Lawrence P. Roberts}
\maketitle

\section{Introduction}

\noindent A knot $K \subset S^{3}$ induces a filtration, $\mathcal{F}_{m}(K)$, on the Heegaard-Floer chain complex
$\widehat{CF}(S^{3})$, \cite{Knot}, whose homology, over the rationals, is $\widehat{HF}(S^{3}) \cong \Q_{\small{(0)}}$. For each $m$ there is an inclusion, $I_{m}$, of chain complexes $\mathcal{F}_{m}(K) \stackrel{I}{\hookrightarrow} \widehat{CF}(S^{3})$. 

\begin{defn}[\cite{OSTa}] Let $K \subset S^{3}$ be a knot, then 
$$
\tau(K) = \min\big\{\,m \in \Z\,\big|\, I_{m,\ast} \mathrm{\ is\ non-trivial}\big\}
$$
\end{defn}
\ \\
\noindent In \cite{OSTa}, \Oz and \Sz proved that 
\ben
\item $\tau(K)$ is an invariant of the concordance class of $K$,
\item $|\tau(K)| \leq g_{4}(K)$, 
\item $\tau(\overline{K}) = -\tau(K)$,
\item $\tau(K \# J) = \tau(K) + \tau(J)$, 
\item $\tau(T_{p,q}) = \disp{\frac{(p-1)(q-1)}{2}}$ where $p, q > 0$ and $T_{p,q}$ is the $(p,q)$-torus knot.
\een 
\ \\
\noindent In this paper, we use four dimensional surgery techniques to provide bounds for the $\tau$-invariant of satellite knots. This approach is prvoided by \cite{OSTa}, Proposition 3.1 wherein properties of surgeries on a knot are shown to compute the $\tau$-invariant up to a bounded error.  The resulting bounds for $\tau$ of a satellite are not particularly strong (they recover the connect sum formula for $\tau$ only up to an error of $\pm 2$, for example), but are quite general. Furthermore, the technique could be useful in other situations. The simplest form of these inequalities is \\

\begin{prop}\label{prop:main}
Let $S_{r}(C,P)$ be the $r$-twisted satellite knot formed from a companion, $C$, in $S^{3}$ and a pattern, $P$, in $S^{1} \times D^{2}$. Let $l$ be the intersection number of $P$ with $D^{2}$, with $P$ oriented so that $l > 0$. Let $n_{+}$ be the minimal number of positive crossings in this intersection number. Define 
$$
D(S_{r}(C,P)) = \tau(S_{r}) - \big( \tau(P) + l\,\tau(C) + \frac{l(l-1)}{2}\,r\big)
$$
then
$$
-n_{+}(P) - l \leq D(S_{r}) \leq n_{+}(P) + l
$$ 
whereas, if $l=0$, we have
$$
-n_{+}(P) - 1 \leq D(S_{r}) \leq n_{+}(P) + 1
$$

\end{prop}

\noindent We note that $P$ will be embedded in $S^{1} \times D^{2}$, with a prescribed framing of $S^{1} \times D^{2}$, making the $r$-twisting well-defined. \\
\ \\
\noindent There are stronger inequalities for certain restricted ranges of $r$.
 
\begin{prop}
Using the notation in the previous proposition, for $r \neq 0$, we have when $l > 0$
$$
\begin{array}{l}
-(1 + l) \leq D(S_{r}) \leq n_{+}(P) + l \hspace{.5in} \mathrm{when\ } r < 2\tau(C) - 1 \\
\ \\
-n_{+}(P) - l \leq D(S_{r}) \leq 1 + l \hspace{.5in} \mathrm{when\ } r > 2\tau(C) + 1 \\
\end{array}
$$
but for $l=0$
$$
\begin{array}{l}
-1 \leq D(S_{r}) \leq n_{+}(P) + 1 \hspace{.5in} \mathrm{when\ } r < 2\tau(C) - 1 \\
\ \\
-n_{+}(P) - 1 \leq D(S_{r}) \leq 1 \hspace{.5in} \mathrm{when\ } r > 2\tau(C) + 1 \\
\end{array}
$$ 
\end{prop} 
\ \\
\noindent Previous results on $\tau(K)$ for satellites revolve around two cases: Whitehead doubles and cables. For cables, M. Hedden
proved the following estimates, based on an analysis of specific Heegaard diagrams

\begin{thm}[Them. 1.2, \cite{Hed1}]
Let $K_{l,lr+1}$ be the $(l,lr+1)$-cable of $K$. Then
$$
l\,\tau(K) + \disp{\frac{lr(l-1)}{2}} \leq \tau(K_{l,lr+1}) \leq l\,\tau(K) + \frac{lr(l-1)}{2} + l - 1
$$
\end{thm}
\ \\
\noindent Furthermore, he gave some cases in which on or other inequality is actually an equality. Furthermore, I. Petkova has used bordered Heegaard Floer homology to compute the $\tau$-invariant explicitly for cables of a kot Floer homologically thin companion, $C$, \cite{Petk}. \\
\ \\
\noindent In \cite{Hed2}, M. Hedden also completely described the $\tau$-invariant for the twisted Whitehead doubles of $K$.  This culminated the work of several authors including C. Livingston and S. Naik, and M. Hedden and P. Ording.\\
\ \\
\noindent Using a different technique, C. Van Cott
also considered cables and discovered, in a slightly different form, that

\begin{thm}[Thrm. 2, \cite{VanC}]
Let $h(n) = \tau(K_{l,n}) - \disp{\frac{(l-1)(n-1)}{2}}$. Then for $n > r$, $n,r$ relatively prime to $l$, 
$$
-(l-1) \leq h(n) - h(r) \leq 0
$$
\end{thm}
\noindent We will make use of a similar result found in \cite{Robe} in Section 4. 

\section{Requisite Heegaard-Floer Results}

\noindent Let $K \subset S^{3}$ be a knot, and let $W_{n}$ be the four dimensional manifold found by attaching a $2$-handle along $K$ with framing $n$. Furthermore, denote the $n$-framed Dehn surgery on $K$ by $S^{3}_{n}(K)$; we will regard $W$ as an oriented cobordism from $S^{3}$ to $S^{3}_{n}(K)$. We can relate $\tau(K)$ to $W$ through\\
\begin{lemma}[Prop 3.1, \cite{OSTa}]
For each $m$, when $n$ is sufficiently large relative to $m$, the cobordism map $\widehat{F}_{n,m}: \widehat{HF}(S^{3}) \ra \widehat{HF}(S^{3}_{-n}(K),[m])$ is non-trivial if $m < \tau(K)$ and is trivial if $m > \tau(K)$, where the map is for the $Spin^{c}$ structure on $W_{-n}$ with $\langle c_{1}(\mathfrak{s}_{m}), \widehat{\Si} \rangle - n = 2\,m$.
\end{lemma}

\noindent Thus the triviality/non-triviality of certain cobordism maps characterizes $\tau(K)$ to within $1$. However, what happens in the lemma when $m = \tau(K)$ is unclear: the cobordism map may or may not be trivial, depending upon $K$. To account for this ambiguity we define

\begin{defn} For a knot $K$, we define the $\tau$-correction, $\mathcal{C}(K)$, to equal $0$ if the the cobordism map in the preceding lemma is non-trivial for $m = \tau(K)$ and to equal $1$ if it is not (assuming $n$ is sufficiently large).
\end{defn}

\noindent Using more recent developments, \cite{IntS}, we can make the previous lemma more precise:

\begin{lemma}\label{prop:chains}
Let $\widehat{F}_{\mathfrak{u}_{i}} : \widehat{HF}(S^{3}) \ra \widehat{HF}(S^{3}_{r}(K),[i])$ be the cobordism map induced from $W_{r}$ and $\langle c_{1}(\mathfrak{s}_{m}), \widehat{\Si} \rangle + r = 2\,i$. For $r \neq 0$, $\widehat{F}_{\mathfrak{u}_{i}}$ is non-trivial if $-\tau(K) + r < i < \tau(K)$. Furthermore, $\widehat{F}_{\mathfrak{u}_{i}} = 0$ if $ i < -\big|\tau(K)\big| + r$ or $i > \big|\tau(K)\big|$ or $r > 0$ and $\tau(K) \leq 0$.
\end{lemma}

\noindent{\bf Proof:} We use the results of \cite{IntS} to compute when $\widehat{F}$ is non-trivial for $n \neq 0$. During the proof of the main theorem in \cite{IntS} (which applies to $HF^{+}$), several computations are performed which allow us to compute the cobordism maps $\widehat{HF}(S^{3}) \ra \widehat{HF}(S^{3}_{r}(K),\mathfrak{s}_{i})$ (see especially sections 4.3, 4.7, and 4.9). This paper tells us that the map $\widehat{F}_{\mathfrak{u}_{i}}$ can be described by an inclusion of chain complexes constructed from the knot Floer homology complex $CFK^{\infty}(K)$. Let $\widehat{C}_{i}$ be the complex

$$
\begindc{\commdiag}[10]
\obj(0,5)[A]{$\disp{\bigoplus_{s \equiv i\mathrm{\,mod\,}|r|}} \widehat{A}_{s}$}
\obj(5,5)[C]{$\ $}
\obj(15,5)[D]{$\ $}
\obj(20,5)[B]{$\disp{\bigoplus_{s \equiv i\mathrm{\,mod\,}|r|}} \widehat{B}_{s}$}
\mor{C}{D}{$\partial_{AB}$} 
\enddc
$$
\ \\
\noindent where $\widehat{A}_{s}$ is the subgroup of $CFK^{\infty}$ determined by $C\{\mathrm{max}(i, j - s) = 0\}$ with differential determined from the quotient complex $A^{+}_{s} = C\{\mathrm{max}(i, j - s) \geq 0\}$, $\widehat{B}_{s}$ is a copy of $C\{i=0\}$ and $\partial_{AB}$ is determined by the maps $\widehat{v}_{s} + [r]\widehat{h}_{s}$ on $\widehat{A}_{s}$. Here $\widehat{v}_{s}: \widehat{A}_{s} \ra \widehat{B}_{s}$ is the map in the sequence
$$
0 \lra C\{i < 0, j =s\} \lra \widehat{A}_{s} \stackrel{\widehat{v}_{s}}{\lra} C\{i = 0\}
$$
and $[r]\widehat{h}_{s}$ is the map $\widehat{h}_{s} = \widehat{D} \circ U^{s} \circ \widehat{H}_{s}: \widehat{A}_{s} \ra \widehat{B}_{s+r}$ where
$$
0 \lra C\{i = 0, j <s\} \lra \widehat{A}_{s} \stackrel{\widehat{H}_{s}}{\lra} C\{j = s\}
$$  
Note that the image of $[r]\widehat{h}_{s}$ lies in $\widehat{B}_{s+r}$, with the shift being determined by $r$. The main result of \cite{IntS} is 
that there is a quasi-isomorphism $\widehat{CF}(S^{3}_{r}(K), [i]) \ra \widehat{C}_{i}$. Under this quasi-isomorphism, the cobordism map $\widehat{F}_{\mathfrak{u}_{i}}$ becomes the inclusion $\widehat{B}_{i} \ra \widehat{C}$, where $\mathfrak{u}_{i}$ is determined by $\langle c_{1}(\mathfrak{u}_{i}), \Si_{K} \rangle + r = 2i$, \cite{IntS}.\\
\ \\
\noindent Assume that $\tau(K) \geq 0$. Then $\widehat{v}_{i} = 0$ when $i < \tau(K)$ since its image is in $\mathcal{F}(K,i)$, and in homology this maps to zero in $C\{i=0\}$. Whereas $\widehat{h}_{i} = 0$ when $i > -\tau(K)$ since it has image in $\mathcal{F}(K,-i)$.  Thus, $H_{\ast}(\widehat{B}_{i}) \cong \Z$ has no non-trivial map into it for $-\tau(K) + r < i < \tau(K)$ and hence the map with $\langle c_{1}(\mathfrak{u}_{i}), \Si_{K} \rangle + r = 2i$ is non-trivial in this set. In particular, $\widehat{F}_{\mathfrak{u}_{i}} \neq 0$ for $i = \tau(K) - 1$ when $r < 2\tau(K) -1$.  At $\tau(K)$ there may be problems, depending upon whether in $\widehat{A}_{\tau(K)}$ there is a horizontal component to the differential applied to the generators mapping non-trivially into $\widehat{HF}(S^{3})$ or not. This is depicted as (where we have taken $\tau(K) = 3$ for concreteness) with a dashed arrow: \\
$$
\begindc{\commdiag}[10]
\obj(4,3)[B1]{$\Z$}
\obj(8,3)[B2]{$\Z$}
\obj(12,3)[B3]{$\Z$}
\obj(16,3)[B4]{$\Z$}
\obj(20,3)[B5]{$\Z$}
\obj(24,3)[B6]{$\Z$}
\obj(28,3)[B7]{$\Z$}
\obj(32,3)[B8]{$\Z$}
\obj(36,3)[B9]{$\Z$}
\obj(40,3)[B10]{$\Z$}
\obj(0,9)[A0]{$A_{-5}$}
\obj(4,9)[A1]{$A_{-4}$}
\obj(8,9)[A2]{$A_{-\tau(K)}$}
\obj(12,9)[A3]{$A_{-2}$}
\obj(16,9)[A4]{$A_{-1}$}
\obj(20,9)[A5]{$A_{0}$}
\obj(24,9)[A6]{$A_{1}$}
\obj(28,9)[A7]{$A_{2}$}
\obj(32,9)[A8]{$A_{\tau(K)}$}
\obj(36,9)[A9]{$A_{4}$}
\obj(40,9)[A10]{$A_{5}$}
\mor{A8}{B8}{$\,$}[1,1]
\mor{A9}{B9}{$\,$}[1,5]
\mor{A10}{B10}{$\,$}[1,5]
\mor{A0}{B3}{$\,$}[1,5]
\mor{A1}{B4}{$\,$}[1,5]
\mor{A2}{B5}{$\,$}[1,1]
\enddc
$$
\ \\
\noindent In fact, this also holds when $\tau(K) < 0$, but the details of the proof are slightly different. Again $\widehat{v}_{i} = 0$ when $i < \tau(K)$ and $\widehat{h}_{i} = 0$ when $i > -\tau(K)$. When $r > 0$, the map will be trivial, since we obtain a complex such as:

$$
\begindc{\commdiag}[10]
\obj(4,3)[B1]{$\Z$}
\obj(8,3)[B2]{$\Z$}
\obj(12,3)[B3]{$\Z$}
\obj(16,3)[B4]{$\Z$}
\obj(20,3)[B5]{$\Z$}
\obj(24,3)[B6]{$\Z$}
\obj(28,3)[B7]{$\Z$}
\obj(32,3)[B8]{$\Z$}
\obj(36,3)[B9]{$\Z$}
\obj(40,3)[B10]{$\Z$}
\obj(0,9)[A0]{$A_{-5}$}
\obj(4,9)[A1]{$A_{-4}$}
\obj(8,9)[A2]{$A_{\tau(K)}$}
\obj(12,9)[A3]{$A_{-2}$}
\obj(16,9)[A4]{$A_{-1}$}
\obj(20,9)[A5]{$A_{0}$}
\obj(24,9)[A6]{$A_{1}$}
\obj(28,9)[A7]{$A_{2}$}
\obj(32,9)[A8]{$A_{-\tau(K)}$}
\obj(36,9)[A9]{$A_{4}$}
\obj(40,9)[A10]{$\ $}
\mor{A2}{B2}{$\,$}[1,1]
\mor{A3}{B3}{$\,$}[1,5]
\mor{A4}{B4}{$\,$}[1,5]
\mor{A5}{B5}{$\,$}[1,5]
\mor{A6}{B6}{$\,$}[1,5]
\mor{A7}{B7}{$\,$}[1,5]
\mor{A8}{B8}{$\,$}[1,5]
\mor{A9}{B9}{$\,$}[1,5]
\mor{A10}{B10}{$\,$}[1,5]
\mor{A8}{B10}{$\,$}[1,1]
\mor{A7}{B9}{$\,$}[1,5]
\mor{A6}{B8}{$\,$}[1,5]
\mor{A5}{B7}{$\,$}[1,5]
\mor{A4}{B6}{$\,$}[1,5]
\mor{A3}{B5}{$\,$}[1,5]
\mor{A2}{B4}{$\,$}[1,5]
\mor{A1}{B3}{$\,$}[1,5]
\mor{A0}{B2}{$\,$}[1,5]
\obj(0,6)[S]{$\,$}
\mor{S}{B1}{$\,$}[1,5]
\enddc
$$
where each of the arrows is onto. Then each of the $\Z$-factors is in the image. For $ i  \leq \tau(K)-1$ or $i > -\tau(K)+r$ this is obvious. For $i=0$, for instance, there is an element in $A_{0}$ mapping onto the generator of the $\Z$ term. If this element does not also map non-trivially to $B_{2}$ then $B_{0}$ is in the image of $A_{0}$. If, however, it does map to $B_{2}$ then there is an element in $A_{2}$ with the same image. We repeat the process with this element. Due to the termination of the slanted arrows at $A_{-\tau(K)+1}$, the process stops. This can be performed for each $B_{i}$ and shows that the inclusion of each $B_{i}$ will be trivial after computing the last stage in the spectral sequence. Of course, the number of arrows involved will depend upon $r$ and $\tau(K)$. \\
\ \\
\noindent However, when $r < 0$ we have a complex such as:

$$
\begindc{\commdiag}[10]
\obj(4,3)[B1]{$\Z$}
\obj(8,3)[B2]{$\Z$}
\obj(12,3)[B3]{$\Z$}
\obj(16,3)[B4]{$\Z$}
\obj(20,3)[B5]{$\Z$}
\obj(24,3)[B6]{$\Z$}
\obj(28,3)[B7]{$\Z$}
\obj(32,3)[B8]{$\Z$}
\obj(36,3)[B9]{$\Z$}
\obj(40,3)[B10]{$\Z$}
\obj(0,9)[A0]{$\,$}
\obj(4,9)[A1]{$\,$}
\obj(8,9)[A2]{$A_{\tau(K)}$}
\obj(12,9)[A3]{$A_{-2}$}
\obj(16,9)[A4]{$A_{-1}$}
\obj(20,9)[A5]{$A_{0}$}
\obj(24,9)[A6]{$A_{1}$}
\obj(28,9)[A7]{$A_{2}$}
\obj(32,9)[A8]{$A_{-\tau(K)}$}
\obj(36,9)[A9]{$A_{4}$}
\obj(40,9)[A10]{$A_{5}$}
\mor{A2}{B2}{$\,$}[1,1]
\mor{A3}{B3}{$\,$}[1,5]
\mor{A4}{B4}{$\,$}[1,5]
\mor{A5}{B5}{$\,$}[1,5]
\mor{A6}{B6}{$\,$}[1,5]
\mor{A7}{B7}{$\,$}[1,5]
\mor{A8}{B8}{$\,$}[1,5]
\mor{A9}{B9}{$\,$}[1,5]
\mor{A10}{B10}{$\,$}[1,5]
\mor{A8}{B6}{$\,$}[1,1]
\mor{A7}{B5}{$\,$}[1,5]
\mor{A6}{B4}{$\,$}[1,5]
\mor{A5}{B3}{$\,$}[1,5]
\mor{A4}{B2}{$\,$}[1,5]
\mor{A3}{B1}{$\,$}[1,5]
\obj(0,3)[B0]{$\,$}
\mor{A2}{B0}{$\,$}[1,5]
\enddc
$$
The situation is no longer clear. When $r < 2\tau(K) - 1$, we obtain a definite non-trivial map with $i = \tau(K) - 1$ since the slanted
arrows will not map into $B_{\tau(K) - 1}$. Pursuing this line of reasoning, we obtain the results in the lemma. $\Diamond$ \\
\ \\
\noindent We removed the case when $\tau(K) <0$ and $\tau(K) < i < -\tau(K)$  from the above lemma due to the possibility of non-trivial maps that can occur if the same element in each $A_{i}$ maps to generators of $B_{i}$ and $B_{i+r}$ and the images of these elements from different $A_{i}$'s cancel sufficiently. For individual knots this will need careful analysis to determine.  \\
\ \\
\noindent We note that in each case in the lemma, the behavior at the ends of the intervals depend upon the knot in question. We will thus have correction terms for each endpoint, as for each framing, but will only need the one at $\tau(K)$ below. We now record some results on these corrections. \\
\ \\
\noindent Let $\mathcal{C}_{r}(K)$ be the correction at $\tau(K)$ for $r$ surgery on $K$. Recall that $\mathcal{C}(K)$ is the correction for sufficiently negative surgeries.

\begin{lemma}
If $r \leq 2\tau(K) - 1$ then $\mathcal{C}_{r}(K) = \mathcal{C}(K)$
\end{lemma}

\noindent{\bf Proof:} We refer to the previous diagrams. First, if $\tau(K) > 0$ then $\mathcal{C}_{r}(K)$ is determined by whether or not $H_{\ast}(A_{\tau(K)}) \stackrel{\widehat{v}}{\lra} \Z_{0}$ is surjective. If it is then $\mathcal{C}_{r}(K) = 1$ and $0$ otherwise. When $r \leq 2\tau(K) - 1$, this map is uninfluenced by the maps
$\widehat{h}$ as they will map into factors to the left of $i = \tau(K)$. (The first non-zero such map is at $-\tau(K)$ and maps to the factor in position $- \tau(K) + r \leq \tau(K) - 1$). Thus for $r \leq 2 \tau(K) - 1$  and $\tau(K) > 0$ we have $\mathcal{C}_{r}(K) = \mathcal{C}(K)$. In fact, if $r \leq 2 \tau(K) - 1$ and $\tau(K) \leq 0$ then $r < 0$ and the map from $H_{\ast}(A_{-\tau(K)})$ to $B_{r - \tau(K)}$ also maps into a factor to the left of $i=\tau(K)$. It again follows, taking into account the possibilities for $\widehat{h}$ that $\mathcal{C}_{r}(K) = \mathcal{C}(K)$.  $\Diamond$\\

\begin{lemma}
If $\tau(K) = g(K)$ then $\mathcal{C}(K) = 1$. If $\tau(K) = -g(K)$ then $\mathcal{C}(K) = 0$. 
\end{lemma}

\noindent {\bf Proof:} The reader should consult the diagrams in the proof of proposition \ref{prop:chains}. If $\tau(K) = g(K)$ then $H_{\ast}(\widehat{A}_{g(K)}) = \Z$, and $\widehat{v}_{g(K)}$ is an isomorphism onto $H_{\ast}(\widehat{B}_{g(K)})$ since $\widehat{A}_{g(K)} \cong \widehat{B}_{g(K)}$ (there cannot be any horizontal component to the differential in $\widehat{A}_{g(K)}$!). Therefore inclusion of $B_{g(K)}$ is trivial in homology, thus the requisite cobordism map is trivial. If $\tau(K) = -g(K)$, let $\xi$ be a class in $\widehat{CF}(K,-g(K))$ which maps isomorphically to $\widehat{HF}(S^{3})$ under the inclusion of chain complexes. If $\partial_{h}\xi = 0$ then $xi' = U^{-g(K)}\xi$ is closed in $C(j=0) \cong \widehat{CF}(S^{3})$. However, in the filtration induced on $\widehat{CF}(S^{3})$ by the $i$-index on $C(j=0)$, which corresponds to the knot $K$ with reversed direction, $\xi'$ occurs in filtration index $g(K)$. No element $\xi' + y$, with $y$ in supported in strictly lower filtration indices, is exact unless $\xi'$ is exact. If $\xi' + y = \partial z$ then $z = \sum z_{i}$. Let $z_{m}$ be the non-zero term supported in the largest index, which is not closed. Since $\xi' + y$ is supported in indices $g(K)$ and below, $z'_{m}$ is closed in $\widehat{CF}(S^{3}, K, m)$, and is not itself exact, hence $\widehat{HF}(K, m) \not\cong \{0\}$. Therefore, $m = g(K)$, and $\xi' = \partial_{0} z'_{m}$ in $\widehat{CF}(K, g(K))$. However, 
$U^{g}\xi' = U^{g}\partial_{0} z'_{m} = \partial_{0} U^{g}z'_{m} = \partial_{0}z$. Then $\xi = \partial_{0} z$ in $\widehat{CF}(K, -g(K))$, which contradicts that $\xi$ generates a summand in $\widehat{HF}(K, -g(K))$, and thus implies that $\partial_{h}(\xi)$ in $\widehat{A}_{-g(K)}$ must be non-zero. Consequently $v_{-g(K), \ast} = 0$. $\Diamond$\\
\ \\
\noindent Types of knots to which this last lemma applies include positive knots, strongly quasi-positive knots, and knots for which a positive 
surgery yields an L-space.\\

\section{Analyzing satellites knots using the lemma}

\noindent Our goal in this section will be to use the previous Heegaard-Floer results to prove the following proposition:

\begin{prop}\label{prop:part2}
Let $S_{r}(C,P)$ be the $r$-twisted satellite knot formed from a companion, $C$, in $S^{3}$ and a pattern, $P$, in $S^{1} \times D^{2}$. Let $l$ be the intersection number of $P$ with $D^{2}$. Let 
$$
D(S_{r}(C,P)) = \tau(S_{r}) - \big( \tau(P) + l\,\tau(C) + \frac{l(l-1)}{2}\,r\big)
$$
then when $r \neq 0$, we have
$$
\begin{array}{l}
-\big(\mathcal{C}(P) + l \,\mathcal{C}(C)\big) \leq D(S_{r}) \leq D(S_{\Delta(r)}(U,P)) + 1 + l\,\mathcal{C}(\overline{C}) \hspace{.5in} \mathrm{when\ } r < 2\tau(C) - 1 \\
\ \\
D(S_{\Delta'(r)}(U,P)) - 1 - l\,\mathcal{C}(C) \leq D(S_{r}) \leq \big(\mathcal{C}(\overline{P}) + l\,\mathcal{C}(\overline{C})\big) \hspace{.5in} \mathrm{when\ } r > 2\tau(C) + 1 \\
\end{array}
$$ 
where $U$ is the unknot, $\Delta(r) = r - 2 \tau(C) - 1 - \mathcal{C}(\overline{C})$, and $\Delta'(r) = r - 2 \tau(C) + 1 + \mathcal{C}(C)$. The first inequality applies to $r=2\tau(C)-1$ when $\mathcal{C}(C) = 0$ and the second set of inequalities applies to $2\tau(C) + 1$ when $\mathcal{C}(\overline{C})= 0$. Furthermore, in all instances:
$$
D(S_{\Delta'(r)}(U,P)) - 1 - l\,\mathcal{C}(C) \leq D(S_{r}) \leq D(S_{\Delta(r)}(U,P)) + 1 + l\,\mathcal{C}(\overline{C})
$$
\end{prop} 

\noindent In the next section we take up the issue of replacing $D(S_{\Delta(r)}(U,P))$ and $D(S_{\Delta'(r)}(U,P))$ by more congenial representations. \\
\ \\
\noindent {\bf Proof:}\\
\ \\
\noindent {\bf I.} Let $C$ be the companion, a knot in $S^{3}$. Let $P$, a knot in $S^{1} \times D^{2}$, be the pattern. Assume that $P$ intersects the oriented $D^{2}$ algebraically a non-negative number of times (re-orient the knot to achieve this). \\
\ \\
\noindent {\bf II.} The $4$-manifold, $W_{r,n}$, found from $r \neq 0$ surgery on $C$ and $-n$ surgery on a geometrically unlinked copy of $P$, can be decomposed as $-n + l^{2} r$ surgery on the $r$-twisted satellite, $S_{r}(C,P)$, union a $2$-handle.\\
\ \\
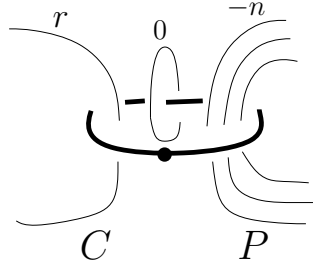
\begin{center}
\includegraphics[scale=0.9]{satellite}
\end{center}
\ \\
\noindent The diffeomorphism can be seen through the diagram above. We may slide $P$ and $C$ over the $0$-framed two handle and then
cancel the one handle with the $0$-framed two handle. $W_{r,n}$ is the cobordism that remains. Alternately, we may slide all the strands of $P$
over $C$, leaving $C$ linking the one handle. Sliding the $0$ handle over $C$ as well, we may cancel $C$. Thus we may first add a handle to the
$r$-twisted satellite, and then add a two handle to obtain the same four manifold. We compute the framing on the satellite:
$$
\left[ \begin{array}{ccc}
			 r & 0 & 1 \\
			 0 & -n & -l \\
			 1 & -l & 0 \\
			 \end{array} \right]
	\lra
\left[ \begin{array}{ccc}
			 r & lr & 1 \\
			 lr & -n + l^{2}r & 0 \\
			 1 & 0 & 0 \\
			 \end{array} \right]
$$
Below, we refer to the four manifold resulting from the surgery on the satellite as $W_{S}$. We will also use the notation $P^{r}$ for a pattern with $r$ twists already in place. We may then find $0$-twisted satellites by the identity
$S_{0}(C,P) = S_{i}(C,P^{-i})$.\\
\ \\
\noindent {\bf III.} The map $\widehat{F}_{W, \mathfrak{s}_{1} \# \mathfrak{s}_{2}}$ is equivalent to $\widehat{F}_{C, \mathfrak{s}_{1}} \otimes \widehat{F}_{P, \mathfrak{s}_{2}}$ under the connect sum isomorphisms.\\
\ \\
\noindent First we describe an isomorphism map, $\widehat{HF}_{\Q}(Y_{1}, \mathfrak{s}_{1}) \otimes \widehat{HF}_{\Q}(Y_{2}, \mathfrak{s}_{2}) \ra \widehat{HF}_{\Q}(Y_{1} \# Y_{2}, \mathfrak{s}_{1} \# \mathfrak{s}_{2})$. Take a pointed Heegaard diagram $(\Si_{1}, {\bf \al}_{1}, {\bf \be}_{1}, z_{1})$ for $Y_{1}$ and a pointed Heegaard diagram, $(\Si_{2}, {\bf \al}_{2}, {\bf \be}_{2}, z_{2})$ for $Y_{2}$, each of which is weakly admissible for the given $Spin^{c}$ structure. The map, on chain complexes, is described as a composition of a map
$$
\widehat{CF}(Y_{1}, \mathfrak{s}_{1}) \otimes \widehat{CF}(Y_{2}, \mathfrak{s}_{2}) \ra \widehat{CF}(\Si, {\bf\al}_{1}{\bf\al}_{2}, {\bf\be}_{1}{\bf\al}_{2}') \otimes \widehat{CF}(\Si, {\bf\be}_{1}{\bf\al}_{2}, {\bf\be}_{1}'{\bf\be}_{2})
$$
where $\Si$ is the connect sum of $\Si_{1}$ and $\Si_{2}$ at $z_{1}, z_{2}$, pointed by $z$ chosen in the connect sum neck, and $(\Si, {\bf \al}_{1}{\bf\al}_{2}, {\bf\be}_{1}{\bf\al}_{2}') \cong Y_{1} \#^{g_{2}} S^{1} \times S^{2}$, $\widehat{CF}(\Si, {\bf\be}_{1}{\bf\al}_{2}, {\bf\be}_{1}'{\bf\be}_{2}) \cong \#^{g_{1}} S^{1} \times S^{2} \# Y_{2}$. On generators, the map takes $${\bf x} \otimes {\bf y} \ra \big({\bf x} \otimes \Theta^{+}_{1} \big) \otimes \big( \Theta^{+}_{2} \otimes {\bf y} \big)$$ where $\Theta^{+}$  is a closed generator for the $\wedge H_{1}$-module structure on $\widehat{HF}(\#^{k} S^{1} \times S^{2}, \mathfrak{s}_{0})$. This map is then composed with the holomorphic triangle map, $\widehat{F}_{\al\be\ga}$, determined by the triple $(\Si, {\bf \al}_{1}{\bf \al}_{2}, {\bf \be}_{1}{\bf \al}_{2}', {\bf \be}_{1}'{\bf \be}_{2}, z)$
and the $Spin^{c}$ structure $\mathfrak{s}_{1} \# \mathfrak{s}_{2}$. That this composition is an isomorphism of chain complexes follows from
the existence of unique ``small'' triangles in the triple diagram, when we choose sufficiently small Hamiltonian isotopes of the attaching circles, indicated by an apostrophe in the data above. For more details, see Section 6 of \cite{3Man}, or Prop. 4.4 of \cite{AbsG}. Over $\Q$, this isomorphism extends to the homologies. We call this map $\widehat{F}_{Y_{1} \# Y_{2}, \mathfrak{s}_{1} \# \mathfrak{s}_{2}}$\\
\ \\
\noindent Using this isomorphism, and the associativity of triple maps, 

\begin{lemma}[Proposition 4.4 of \cite{AbsG}]
The map $\widehat{F}_{Y \# Z, \mathfrak{s} \# \mathfrak{t}}$ is independent of the Heegaard diagrams used for $Y$ and $Z$. Moreover, if $W$ is cobordism from $Y$ to $Y'$, equipped with a $spin^{c}$ structure $\mathfrak{u}$, restricting to the ends appropriately, then the following diagram is commutative:
$$
\begindc{\commdiag}[10]
\obj(0,0)[I00]{$\widehat{HF}(Y', \mathfrak{s}') \otimes \widehat{HF}(Z, \mathfrak{t})$}
\obj(16, 0)[I10]{$\widehat{HF}(Y' \# Z, \mathfrak{s}'\#\mathfrak{t})$}
\obj(0,7)[I01]{$\widehat{HF}(Y, \mathfrak{s}) \otimes \widehat{HF}(Z, \mathfrak{t})$}
\obj(16,7)[I11]{$\widehat{HF}(Y \# Z, \mathfrak{s}\#\mathfrak{t})$}
\mor{I00}{I10}{$\widehat{F}_{Y' \# Z, \mathfrak{s}' \# \mathfrak{t}}$}
\mor{I11}{I10}{$\widehat{F}_{W \# (Z \times I), \mathfrak{u}\#\mathfrak{t}}$}
\mor{I01}{I00}{$\widehat{F}_{W,\mathfrak{u}} \otimes \mathrm{Id}$}
\mor{I01}{I11}{$\widehat{F}_{Y \# Z, \mathfrak{s} \# \mathfrak{t}}$}
\enddc
$$
\end{lemma}

\noindent This is simply altering Prop. 4.4 to apply to $\widehat{HF}$. Obviously, the order of the factors does not matter, and we can apply the lemma also to a non-trivial cobordism $Z \ra Z'$. To complete step {\bf III} we now apply this lemma twice: once to $S^{3} \# S^{3} \ra S^{3}_{r}(C) \# S^{3}$ and the second time to $S^{3}_{r}(C) \# S^{3} \ra S^{3}_{r}(C) \# S^{3}_{-n}(P)$. The cobordism maps on the right compose on homology to give the overall cobordism map (since $S^{3}_{r}(C)$ is a rational homology sphere) for the $spin^{c}$ structure determined as above. On the left, we obtain the map $\widehat{F}_{C, \mathfrak{s}_{1}} \otimes \widehat{F}_{P, \mathfrak{s}_{2}}$. \\
\ \\
\noindent {\bf IV.} On the other hand $\partial W_{S}$ is a rational homology sphere for sufficiently large $n$. Thus, we can decompose
$\widehat{F}_{W} = \widehat{F}_{h} \circ \widehat{F}_{S}$. The $Spin^{c}$ structures also decompose, and are determined by how they restrict to 
the two parts. We thus have that $\widehat{F}_{W} \neq 0$ implies that $\widehat{F}_{S}\neq 0$. We now undertake to determine how this
occurs relative to the $Spin^{c}$ structures. \\
\ \\
\noindent {\bf V.} We have the freedom above to vary $n$ as we wish. Assume that $\widehat{F}_{C} \neq 0$ for $\langle c_{1}(\mathfrak{s}_{1}), \widehat{\Si}_{1} \rangle + r = 2\,k$ and that $\langle c_{1}(\mathfrak{s}_{2}), \widehat{\Si} \rangle - n = 2\,m$, for $n$ sufficiently large. Then $-n + l^{2}\,r$ can be made sufficiently negative for the surgery lemma to apply to both $P$ and $S_{r}$. We have that
$$
\langle c_{1}(\mathfrak{s}), \widehat{\Si}_{2} + l\widehat{\Si}_{1} \rangle = (2m + n) + l\,(2k - r) = 2(m + l\,k) + n -l\,r + l^{2}\,r-l^{2}\,r
$$
\noindent Re-organizing yields
$$ 
\langle c_{1}(\mathfrak{s}), \widehat{\Si} \rangle + (-n + l^{2}\,r) = 2(m + l\,k) + l(l-1)\,r
$$
\noindent If $m < \tau(P)$ then we must have that the map for this $Spin^{c}$ structure is non-trivial. Thus
$$
m + l\,k + \frac{l(l-1)}{2}\,r \leq \tau(S_{r})
$$
This will hold for $m = \tau(P) - 1$ (or perhaps $\tau(P)$ depending on the value of $\mathcal{C}(P)$). If we let $\overline{R}$ be the largest value for which $r$ surgery on $C$ has a non-zero map satisfying the relation above, then we have
$$
\tau(P) + l\,\overline{R} + \frac{l(l-1)}{2}\,r \leq \tau(S_{r}) + \mathcal{C}(P)
$$
\ \\
\noindent From the previous sections we know that $\widehat{F}_{r,i}^{C}$ is non-trivial for $-\tau(C) + r < i < \tau(C)$ and sometimes we may be able to extend these to the endpoints. So for $r < 2\tau(C) - 1$, we may take $i = \tau(C) - 1$ and have a non-zero map. Sometimes, we can extend to $i = \tau(C)$, and then we may use $r = 2\tau(C) - 1$ as well. We then have
$$
\tau(P) + l\,\tau(C) + \frac{l(l-1)}{2}\,r \leq \tau(S_{r}) + \mathcal{C}(P) + l\,\mathcal{C}(C) \hspace{.5in} \mathrm{when\ } r < 2\tau(C) - 1 
$$
with the possibility of including $r = 2\tau(C) - 1$ if $\mathcal{C}(C) = \mathcal{C}_{r}(C) = 0$.\\
\ \\
Let 
$$
T(S_{r}) = \tau(P) + l\,\tau(C) + \frac{l(l-1)}{2}\,r
$$
and
$$
D(S_{r}) = \tau(S_{r}) - T(S_{r})
$$
then we have established one-sided bound:
$$
\begin{array}{l}
-\big(\mathcal{C}(P) + l \,\mathcal{C}(C)\big) \leq D(S_{r}) \hspace{.5in} \mathrm{when\ } r < 2\tau(C) - 1 \\
\end{array}
$$
including the endpoint for $r$ if $\mathcal{C}(C) = 0$. \\
\ \\
\noindent {\bf VI.} The crucial observation in finding bounds for the other side is 
$$
\tau(S_{m}(\overline{C}, S_{r}(C,P))) = \tau(S_{r + m}(U, P)) = \tau(P^{r+m})
$$
This follows from the observation that \\

\begin{lemma}\label{lem:satcon}
The knot $S_{m}(\overline{C}, S_{r}(C,P))$ is concordant to $S_{m + r}(U, P)$.
\end{lemma}

\noindent {\bf Proof:} Since $C \# \overline{C}$ is a ribbon knot, we can find a slice disk for it. If we trivialize a neighborhood of this disc to obtain a region of $B^{4}$ diffeomorphic to $D^{2} \times D^{2}$. If we take $D^{2} \times \{0\}$ to be the slice disc, we can construct $m+k$ parallel copies by choosing $m+r$ points, $x_{i}$, in the second factor, and taking the image of $D^{2} \times \{x_{i}\}$ under the diffeomorphism for each $i$. In $S^{3}$ this gives a link formed by $m+r$ parallel copies of $C \# \overline{C}$. Each copy is a longitude since it bounds a disc disjoint from the slice disc. We place this configuration closed to $S_{m+r}(U,P)$, and orient the longitudes in such a way that we can perform $n$ band sums and obtain an oriented knot. This knot is the same as $S_{m}(\overline{C}, S_{r}(C,P))$. $\Diamond$ \\
\ \\
\noindent Then we immediately have that
$$
-\big(\mathcal{C}(S_{r}) + l\,\mathcal{C}(\overline{C})\big) \leq \tau(S_{m}(\overline{C}, S_{r}(C,P))) - T(S_{m}(\overline{C}, S_{r}(C,P))) \hspace{.5in} \mathrm{when\ } m < 2\tau(\overline{C}) - 1 \\
$$
but
$$
T(S_{m}(\overline{C}, S_{r}(C,P))) = \tau(S_{r}) + l\,\tau(\overline{C}) + \frac{l(l-1)}{2}\,m  
$$
therefore
$$
\tau(S_{r}) \leq \tau(P^{r + m}) -  \big(l\,\tau(\overline{C}) + \frac{l(l-1)}{2}\,m\big) + \mathcal{C}(S_{r}) + l\,\mathcal{C}(\overline{C}) \hspace{.5in} \mathrm{when\ } m < -2\tau(C) - 1 \\
$$
As above, if $\mathcal{C}(\overline{C}) = 0$ we can extend to $m = -2\tau(C) -1$. We choose $m = - 2\tau(C) - 1 - \mathcal{C}(\overline{C})$. Rearranging and simplifying we obtain:
$$
\tau(S_{r}) \leq \tau(P^{r - 2 \tau(C) - 1 - \mathcal{C}(\overline{C})}) +  l\,\tau(C) + \frac{l(l-1)}{2}\,(2\tau(C) + 1 + \mathcal{C}(\overline{C})) + \mathcal{C}(S_{r}) + l\,\mathcal{C}(\overline{C}) \\
$$
Let $\Delta(r) = r - 2 \tau(C) - 1 - \mathcal{C}(\overline{C})$, then 
$$
\tau(S_{r}) - T(S_{r}) \leq \tau(P^{\Delta(r)}) - \tau(P) - \frac{l(l-1)}{2}\,\Delta(r) + \mathcal{C}(S_{r}) + l\,\mathcal{C}(\overline{C})
$$
We thus have
$$
\tau(S_{r}) - T(S_{r}) \leq \tau(P^{\Delta(r)}) - T(P^{\Delta(r)}) + 1 + l\,\mathcal{C}(\overline{C}) 
$$
where we have replaced $\mathcal{C}(S_{r})$ with $1$ to ensure each of the inequalities holds regardless of information about $S_{r}$. Note that
this inequality always applies, regardless of the value of $r$, since it depends only upon the choice of $m$.\\
\ \\
\noindent  Taken with the previous inequality, the above yields:
$$
-\big(\mathcal{C}(P) + l \,\mathcal{C}(C)\big) \leq D(S_{r}) \leq D(P^{\Delta(r)}) + 1 + l\,\mathcal{C}(\overline{C}) \hspace{.5in} \mathrm{when\ } r < 2\tau(C) - 1 \\
$$
which also applies at $r=2\tau(C) - 1$ if $\mathcal{C}(C) = 0$.\\
\ \\
\noindent {\bf VII.} For $r > 2\tau(C) + 1$, with the same caveats about corrections, consider the
satellite knot $S_{-r}(\overline{C}, \overline{P})$ which is the mirror of $S_{r}(C, P)$. Then $\tau(S_{-r}(\overline{C}, \overline{P})) = - \tau(S_{r})$ and $T(S_{-r}(\overline{C}, \overline{P}) = - T(S_{r}(C,P))$. Thus $D(\overline{S_{r}}) = - D(S_{r})$ for any satellite. In addition $\Delta_{C}(r) = -\Delta'_{\overline{C}}(-r)$ where $\Delta'_{C}(r) = r - 2 \tau(C) + 1 + \mathcal{C}(C)$. Since $-r \leq 2\tau(\overline{C}) - 1$ means $r \geq 2\tau(C) + 1$, we see that when applying the preceding inequality to  $S_{-r}(\overline{C}, \overline{P})$ when $r > 2\tau(C) + 1$ we
have
$$
-\big(\mathcal{C}(\overline{P}) + l \,\mathcal{C}(\overline{C})\big) \leq D(\overline{S_{r}}) \leq D(\overline{P}^{\Delta_{\overline{C}}(-r)}) + 1 + l\,\mathcal{C}(C) \hspace{.5in} \mathrm{when\ } r > 2\tau(C) + 1 \\
$$
or since $\overline{P}^{\Delta_{\overline{C}}(-r)} = \overline{P^{\Delta'_{C}(r)}}$ and using the change in sign for $D(S_{r})$ we have
$$
D(P^{\Delta'(r)}) - 1 - l\,\mathcal{C}(C) \leq D(S_{r}) \leq \big(\mathcal{C}(\overline{P}) + l\,\mathcal{C}(\overline{C})\big) \hspace{.5in} \mathrm{when\ } r > 2\tau(C) + 1 \\
$$
We can include $r= 2\tau(C) + 1$ if $\mathcal{C}(\overline{C}) = 0$. Here the left hand inequality holds regardless of the value of $r$. This 
concludes the proof of the proposition. $\Diamond$

\section{Tidying up the inequalities}

\noindent We now wish to clean up the results of the previous section. In particular, we would like to compute $D(S_{\Delta(r)}(U,P))$ in some 
simpler manner. The key will be the following proposition, found in \cite{Robe}, whose proof mimics Van Cott's arguments in \cite{VanC}. 

\begin{prop}
Let the orientation on $P$ be such that $l \geq 0$ and let
$$
g(r) = \tau(S_{r}(C,P)) - \frac{l(l-1)}{2}r
$$
Let $n_{+}$, and $n_{-}$ be the number of strands of $P$ intersecting the
oriented copy of $D^{2}$ positively and negatively, respectively. Then if $s > r$ and $n_{+} > n_{-}$
$$
-(n_{+} - 1) \leq g(s) - g(r) \leq n_{-}
$$
while when $n_{+} = n_{-}$ we have
$$
-n_{+} \leq g(s) - g(r) \leq (n_{-} - 1)
$$
\end{prop}

\noindent With this proposition in hand, we can complete the proof of Proposition \ref{prop:main}:

\begin{lemma} \label{prop:part} When $n_{+} > n_{-}$
$$
-(n_{+}(P) - 1) \leq D(S_{\Delta'(r)}(U,P)) \leq n_{-}(P) \hspace{.5in} r \geq 2\tau(C) + 1\\
$$
$$
$$
$$
-n_{-}(P) \leq D(S_{\Delta(r)}(U,P) \leq n_{+}(P) - 1 \hspace{.5in} r \leq 2\tau(C) - 1
$$
while when $n_{+} = n_{-}$ we have
$$
- n_{+}(P) \leq D(S_{\Delta'(r)}(U,P)) \leq (n_{-}(P) - 1) \hspace{.5in} r \geq 2\tau(C) + 1\\
$$
$$
$$
$$
-(n_{-}(P) - 1) \leq D(S_{\Delta(r)}(U,P) \leq n_{+}(P) \hspace{.5in} r \leq 2\tau(C) - 1
$$
\end{lemma}

\noindent{\bf Proof:} Since $C$ is the unknot we have
$
D(S_{t}(U,P)) = \tau(S_{t}(U,P)) - \tau(P) - \frac{l(l-1)}{2}t
$
which in turn becomes $g(t) - g(0)$, using the notation in the previous proposition. Now $\Delta'(r) \geq 0$ since we only use it
when $r \geq 2 \tau(C) - 1$. Thus,  $-(n_{+} - 1) \leq g(\Delta'(r)) - g(0) \leq n_{-}$ when $r > 2 \tau(C) + 1$ and $n_{+} > n_{-}$ ($- n_{+} \leq g(\Delta'(r)) - g(0) \leq (n_{-} - 1)$ when $n_{+} = n_{-}$). Furthermore, when $r \leq 2 \tau(C) + 1$
we have that $\Delta(r) \leq 0$, but $-(n_{+} - 1) \leq g(0) - g(\Delta(r)) \leq n_{-}$ in this case (when $n_{+} = n_{-}$ this becomes $ -n_{+} \leq g(0) - g(\Delta(r)) \leq (n_{-} - 1)$). Multiplying by $-1$, we can reverse the inequalities. The case when $l=0$ is identical, except that the bounds change as in the previous proposition. $\Diamond$.\\
\ \\
\noindent The inequality in the previous section then becomes:
\begin{prop}
When $r \neq 0, l > 0$, we have
$$
\begin{array}{l}
-\big(\mathcal{C}(P) + l \,\mathcal{C}(C)\big) \leq D(S_{r}) \leq n_{+}(P) + l\,\mathcal{C}(\overline{C}) \hspace{.5in} \mathrm{when\ } r < 2\tau(C) - 1 \\
\ \\
-n_{+}(P) - l\,\mathcal{C}(C) \leq D(S_{r}) \leq \big(\mathcal{C}(\overline{P}) + l\,\mathcal{C}(\overline{C})\big) \hspace{.5in} \mathrm{when\ } r > 2\tau(C) + 1 \\
\end{array}
$$ 
If $\mathcal{C}(C) = 0$ then the first inequality also applies for $r=2\tau(C) - 1$, while if $\mathcal{C}(\overline{C}) = 0$ then
the second inequality applies at $r= 2\tau(C) + 1$. Furthermore, for all $r$ we have
$$
-n_{+}(P) - l\,\mathcal{C}(C) \leq D(S_{r}) \leq n_{+}(P) + l\,\mathcal{C}(\overline{C})
$$
\end{prop} 
 \ \\
\noindent {\bf Proof:} We substitute one side of the inequalities from Proposition \ref{prop:part} into the inequalities in Proposition \ref{prop:part2}. All that remains is the inequalities that hold in general. We know that, for all $r$, 
$$
D(S_{r}) \leq D(S_{\Delta(r)}(U,P)) + 1 + l\,\mathcal{C}(\overline{C})
$$
and $D(S_{\Delta(r)}(U,P)) + 1 + l\,\mathcal{C}(\overline{C}) \leq n_{+}(P) + l\,\mathcal{C}(\overline{C})$ for $r \leq 2 \tau(C) + 1$. But for $r > 2\tau(C) + 1$ we also have $D(S_{r}) \leq \big(\mathcal{C}(\overline{P}) + l\,\mathcal{C}(\overline{C})\big)$. Since $n_{+}(P) \geq 1 \geq \mathcal{C}(\overline{P})$, the inequality on the right holds for all $r$. A similar argument establishes the result for the inequality on the left. $\Diamond$.
\ \\
\noindent To obtain the propositions in the introduction, we set $\mathcal{C}(C) = 1$, which is the worst case for both sides of the inequalities
above. Finally, we address the case when $l=0$. This is identical to that above, but with the different bounds we obtain
\noindent The inequality in the previous section then becomes:
\begin{prop}
When $r \neq 0, l = 0$, we have
$$
\begin{array}{l}
-\mathcal{C}(P) \leq D(S_{r}) \leq n_{+}(P) + 1 \hspace{.5in} \mathrm{when\ } r < 2\tau(C) - 1 \\
\ \\
-n_{+}(P) - 1 \leq D(S_{r}) \leq \mathcal{C}(\overline{P}) \hspace{.5in} \mathrm{when\ } r > 2\tau(C) + 1 \\
\end{array}
$$ 
If $\mathcal{C}(C) = 0$ then the first inequality also applies for $r=2\tau(C) - 1$, while if $\mathcal{C}(\overline{C}) = 0$ then
the second inequality applies at $r= 2\tau(C) + 1$. Furthermore, for all $r$ we have
$$
-n_{+}(P) - 1 \leq D(S_{r}) \leq n_{+}(P) + 1
$$
\end{prop} 
 
\section{Special Cases}

\noindent Below, we assume that $r \neq 0$.

\subsection{When $l=0$ and $P$ is an unknot, when considered in $S^{3}$:} A calculation shows that $\mathcal{C}(P) = 0$ in this case. So we obtain
$$
\begin{array}{l}
0 \leq D(S_{r}) \leq n_{+}(P) + 1 \hspace{.5in} \mathrm{when\ } r < 2\tau(C) - 1 \\
\ \\
-n_{+}(P) - 1 \leq D(S_{r}) \leq 0 \hspace{.5in} \mathrm{when\ } r > 2\tau(C) + 1 \\
\end{array}
$$ 
which conforms to the behavior found for Whitehead doubles in \cite{Hed2}, \cite{LivN}.

\subsection{When $l=1$:} We obtain the inequalities:

$$
\begin{array}{l}
-\big(\mathcal{C}(P) + \mathcal{C}(C)\big) \leq \tau(S_{r}) - \tau(P) - \tau(C) \leq  n_{+}(P) + \mathcal{C}(\overline{C}) \hspace{.5in} \mathrm{when\ } r < 2\tau(C) - 1 \\
\ \\
- n_{+}(P) - \mathcal{C}(C) \leq \tau(S_{r}) - \tau(P) - \tau(C) \leq \big(\mathcal{C}(\overline{P}) + \mathcal{C}(\overline{C})\big) \hspace{.5in} 
\mathrm{when\ } r > 2\tau(C) + 1 \\
\end{array}
$$
\ \\
\noindent If $l=1$ both algebraically and geometrically, then $S_{r} \cong P \# C$ for all $r$. These inequalities almost give the additivity formula under connect sum -- since $n_{+}(P) = 1$ -- but not quite. With some effort, we could replace the correction terms, or simply replace them by $1$'s. In the latter case, we obtain  

$$
\begin{array}{l}
- 2 \leq \tau(S_{r}) - \tau(P) - \tau(C) \leq  n_{+}(P) + 1 \hspace{.5in} \mathrm{when\ } r < 2\tau(C) - 1 \\
\ \\
- (n_{+}(P) + 1) \leq \tau(S_{r}) - \tau(P) - \tau(C) \leq 2 \hspace{.5in} 
\mathrm{when\ } r > 2\tau(C) + 1 \\
\end{array}
$$
\ \\
\noindent These inequalities apply independently of $r$.
\ \\
\subsection{When $P$ is a specific unknot:}

\noindent Let $P$ be the closure of the braid $\sigma_{1}\sigma_{2}\ldots\sigma_{l-1}$. Then $\tau(P) = 0$, $\mathcal{C}(P) = 0$, and $n_{+}(P) = l$. Consequently, we have the inequalities:

$$
\begin{array}{l}
l\,\tau(C) + \disp{\frac{l(l-1)}{2}\,r} -l \leq \tau(S_{r}) \leq  l\,\tau(C) + \disp{\frac{l(l-1)}{2}\,r} + 2l \hspace{.5in} \mathrm{when\ }  r < 2\tau(C) - 1 \\
\ \\
l\,\tau(C) + \disp{\frac{l(l-1)}{2}\,r} -2l \leq \tau(S_{r}) \leq  l\,\tau(C) + \disp{\frac{l(l-1)}{2}\,r} + l \hspace{.5in} \mathrm{when\ } r > 2\tau(C) + 1 \\
\end{array}
$$
\ \\
These are similar to those in \cite{Hed1}.\\
\ \\
\subsection{When $C$ is the unknot:} We write $P^{m} = S_{m}(U,P)$. This is just shorthand for adding full twists to
a collection of parallel strands in $P$. Then 

$$
\begin{array}{l}
-\mathcal{C}(P) \leq \tau(P^{r}) - \tau(P) - \disp{\frac{l(l-1)}{2}\,r} \leq n_{+}(P) \hspace{1in} \mathrm{when\ } r < - 1 \\
\ \\
-n_{+}(P) \leq \tau(P^{r}) - \tau(P) -  \disp{\frac{l(l-1)}{2}\,r} \leq \mathcal{C}(\overline{P})  \hspace{1in} \mathrm{when\ } r >  + 1 \\
\end{array}
$$ 
\ \\
\noindent If $P$ is the closure of a $l$ stranded braid then
$n_{+}(P) = l$ and we obtain

$$
\begin{array}{l}
\disp{\frac{l(l-1)}{2}\,r} - 1 \leq \tau(P^{r}) - \tau(P) \leq \disp{\frac{l(l-1)}{2}\,r} + l \hspace{1in} \mathrm{when\ } r < - 1 \\
\ \\
\disp{\frac{l(l-1)}{2}\,r} - l \leq \tau(P^{r}) - \tau(P) \leq \disp{\frac{l(l-1)}{2}\,r} + 1  \hspace{1in} \mathrm{when\ } r >  + 1 \\
\end{array}
$$
\ \\
\noindent These are similar to the results in section 4 of \cite{VanC}.

\end{document}